\documentclass[pre,10pt,twocolumn]{revtex4}%
\usepackage{amsfonts}
\usepackage{amsmath}
\usepackage{amssymb}
\usepackage{graphicx}%
\newtheorem{theorem}{Theorem}

\newtheorem{definition}[theorem]{Definition}

\newenvironment{proof}[1][Proof]{\noindent\textbf{#1.} }{\ \rule{0.5em}{0.5em}}
\begin{document}
\title[Combining tensor products of graphs]{Tensor 2-sums and entanglement}
\author{Sandi Klav\v{z}ar}
\affiliation{Faculty of Mathematics and Physics, University of Ljubljana, Jadranska 19,
1000 Ljubljana, Slovenia; Faculty of Natural Sciences and Mathematics,
University of Maribor, Koro\v{s}ka 160, 2000 Maribor, Slovenia}
\author{Simone Severini}
\affiliation{Department of Physics and Astronomy, University College London, Gower Street,
WC1E 6BT London, United Kingdom}

\begin{abstract}
To define a minimal mathematical framework for isolating some of the
characteristic properties of quantum entanglement, we introduce a
generalization of the tensor product of graphs. Inspired by the notion of a
density matrix, the generalization is a simple one: every graph can be
obtained by addition modulo two, possibly with many summands, of tensor
products of adjacency matrices. In this picture, we are still able to prove a
combinatorial analogue of the Peres-Horodecki criterion for testing separability.

\end{abstract}
\maketitle

\section{Introduction}

In this paper we attempt to define a minimal mathematical framework for
isolating some of the characteristic properties of quantum entanglement. The
proposed model is such that we hope the paper is amenable to be read by two
audiences with different interests: the audience interested in algebraic and
structural graph theory and the audience interested in entanglement theory.

The tensor product has a fundamental role in the standard \ formulation of
quantum mechanics as the axiomatically designed operation for combining
Hilbert spaces associated to the parties of a distributed quantum mechanical
system (see, \emph{e.g.}, \cite{di}). The definition of entanglement is so
essentially dependent on the tensor product, in a way that we can speak about
entanglement only in the presence of this operation. In the light of such a
fact, mathematical criteria for detecting and classifying entanglement are
mainly based on tools that give information, in most of the cases only
partial, about the tensor product structure of quantum states or their
dynamical operators.

It is plausible that some characteristic properties of significance in the
quantum context remain associated to the tensor product even when we
impoverish the mathematical structure used in quantum mechanics itself. In
different terms, it is conceivable that certain properties of entanglement can
be studied outside quantum mechanics, in a more controlled mathematical
laboratory, where we keep features designated as essential and throw away
redundant or \textquotedblleft less important\textquotedblright\ material. It
is clear that such an experiment would imply a loss of some kind.

The goal of this note is to define a toy-setting with \textquotedblleft fake
quantum states\textquotedblright, which are still combined by using the notion
of tensor product. We do not ask whether we can actually define a physical
theory with a state-space equivalent to the one proper of quantum mechanics,
but obtained with a restricted mathematical tool-box. As we have stated above,
what we do aim for is to picture a scenery with mathematical objects poorer
than general quantum mechanical states, but still exhibiting some of their
characteristic features.

The idea is then to distill a likely analogue of entanglement but in a slimmer
mathematical setting. Specifically, we should be able to: \emph{(i)} define an
operation for mixing states, that is, to obtain statistical mixtures of pure
states; and to \emph{(ii)} define an operation for combining states. Labeled
graphs provide a versatile language for this intent: we mix by sum modulo two
of adjacency matrices; we combine by tensor product of graphs. The latter
operation is well-studied in graph theory. Indeed, it appeared in many
different contexts with a number of equivalent names: tensor product in
\cite{abha-08,zhzh-07}, but also direct product \cite{brsp-08,ha-08} and
categorical product \cite{ra-05,ta-05}, just to mention the most important
ones. See also the recent papers~\cite{ha-09,impi-08,mavu-08}, while for a
general treatment of this graph product we refer to the book~\cite{imrich00}.

Graph tensor products have found a variety of applications. For example, let
us just mention here that recently, Leskovec \emph{et al.}~\cite{les} proposed
tensor powers of graphs for modeling complex networks. The Kronecker product
not only allows an investigation using analytical tools (which is not
surprising since this is a is well understood operation), but the construction
itself results very close to real-world networks.

The remainder of this paper is organized as follows. In the next section we
provide the required preliminary definitions . Then, in Section~\ref{ch},
Theorem~\ref{thm:2char} gives a combinatorial characterization of tensor
2-sums. In Section~\ref{tra}, Theorem~\ref{thm:transpose} gives a
combinatorial analogue of the Peres-Horodecki criterion (see, \emph{e.g.},
\cite{nc}) for testing separability. The concluding section contains several
topics for further research and related problems. In particular, it is an open
question to establish computational complexity results concerning the
recognition problem of tensor 2-sums.

\section{Definitions\label{def}}

We consider graphs with a finite number of vertices, without multiple edges,
and without self-loops. The tensor product of graphs (see the figure for two
examples) is defined as follows:

\begin{definition}
The \emph{tensor product}, $K=G\otimes H$, of graphs $G=(V(G),E(G))$ and
$H=(V(H),E(H))$ is the graph with vertex-set $V(K)=V(G)\times V(H)$ and
$\{\left(  g,h\right)  ,(g^{\prime},h^{\prime})\}\in E(K)$ if and only if
$\left\{  g,g^{\prime}\right\}  \in E(G)$ and $\left\{  h,h^{\prime}\right\}
\in E(H)$.
\end{definition}

%

%TCIMACRO{\FRAME{fhFU}{3.4662in}{1.5281in}{0pt}{\Qcb{The graph on the left is
%the tensor product of two complete graphs on four and three vertices,
%respectively. The graph in the middle (consisting of two connected components)
%is the tensor product of a cycle on four vertices and a path on three
%vertices. The graph on the right is their 2-sum.}}{}{figura1.eps}%
%{\special{ language "Scientific Word";  type "GRAPHIC";
%maintain-aspect-ratio TRUE;  display "USEDEF";  valid_file "F";
%width 3.4662in;  height 1.5281in;  depth 0pt;  original-width 6.8753in;
%original-height 3in;  cropleft "0";  croptop "1";  cropright "1";
%cropbottom "0";  filename 'figura1.eps';file-properties "XNPEU";}}}%
%BeginExpansion
\begin{figure}
[h]
\begin{center}
\includegraphics[
height=1.5281in,
width=3.4662in
]%
{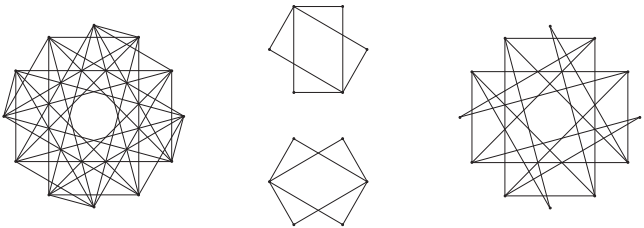}%
\caption{The graph on the left is the tensor product of two complete graphs on
four and three vertices, respectively. The graph in the middle (consisting of
two connected components) is the tensor product of a cycle on four vertices
and a path on three vertices. The graph on the right is their 2-sum.}%
\end{center}
\end{figure}
%EndExpansion

Note that the product graph $K$ is undirected, since $\{\left(  g,h\right)
,(g^{\prime},h^{\prime})\}\in E(K)$ if and only if we have $\{(g^{\prime
},h^{\prime}),(g,h)\}\in E(K)$. Let $G$ be a graph with $V(G)=\{g_{1}%
,g_{2},\ldots,g_{n}\}$. Recall that the \emph{adjacency matrix} $A(G)$ of $G$
is an $n\times n$ matrix with $A(G)_{i,j}=1$ if $\{g_{i},g_{j}\}\in E(G)$ and
$A(G)_{i,j}=0$, otherwise. Note that $A(G)$ is symmetric and that its labeling
depends on the ordering of the vertices of $G $. Let $H$ be another graph with
$V(H)=\{h_{1},h_{2},\ldots,h_{m}\}$. Then, unless stated otherwise, the
adjacency matrix of the tensor product $G\otimes H$ will be understood with
respect to the lexicographic ordering of $V(G\otimes H)$: $(g_{1}%
,h_{1}),\ldots,(g_{1},h_{m}),(g_{2},h_{1}),\ldots,(g_{2},h_{m}),\ldots
,(g_{n},h_{m})$. Under this agreement, the following statement is a well-known
useful fact: if $K=G\otimes H$ then $A(K)=A(G)\otimes A(H)$.

Our generalization of the tensor product of graphs requires an additional
operation that reminds of the symmetric difference, but producing a graph on
the same vertex set of the operands:

\begin{definition}
Let $G$ and $H$ be graphs with $V(G)=V(H)$. The \emph{sum modulo }$2$ or
\emph{2-sum} for short, $K=G\oplus H$, of $G$ and $H$ is the graph with the
same vertex set as $G$ (and as $H$) such that $\{u,v\}\in E(K)$ if and only if
either \emph{(i) }$\{u,v\}\in E(G)$ and $\{u,v\}\notin E(H)$ or \emph{(ii)
}$\{u,v\}\notin E(G)$ and $\{u,v\}\in E(H)$.
\end{definition}

The right-hand side of figure gives an example of the 2-sum of two graphs (in
fact, of two tensor product graphs). The graph $K=G\oplus H$ has adjacency
matrix with $ij$-th entry $\left(  A(K)\right)  _{i,j}=((A(G))_{i,j}%
+(A(H))_{i,j})\operatorname{mod}2$.

We are now ready to give the following definition, where by a nontrivial graph
we mean a graph with a least one edge:

\begin{definition}
A graph $K$ is a \emph{tensor 2-sum} if there exists a positive integer $l$,
nontrivial graphs $G_{1},...,G_{l}$, and nontrivial graphs $H_{1},...,H_{l}$,
such that
\[
K=\bigoplus_{k=1}^{l}(G_{k}\otimes H_{k})\,.
\]

\end{definition}

Notice that the case $l=1$ reduces to the standard tensor product.

\section{Characterization\label{ch}}

Let $\mathcal{K}(p,q)$ be the set of graphs that are a tensor 2-sum in which
the factors of the corresponding tensor products are of size $p$ and $q$,
respectively. Hence $|V(K)|=n=pq$ for $K\in\mathcal{K}(p,q)$. Notice that if
$K\in\mathcal{K}(p,q)$ then
\[
|E(K)|\leq\binom{pq}{2}-q\binom{p}{2}-p\binom{q}{2}=2\binom{p}{2}\binom{q}%
{2},
\]
with the equality if and only if $K=K_{p}\otimes K_{q}$. Let $p\geq2$ and
$q\geq2$ be arbitrary but fixed integers. Let $G$ and $H$ be arbitrary graphs
on $p$ and $q$ vertices, respectively. For our purposes, we may assume that
$V(G)=\{g_{1},\ldots,g_{p}\}$ for an arbitrary graph $G$ on $p$ vertices and
$V(H)=\{h_{1},\ldots,h_{q}\}$ for an arbitrary graph $H$ on $q$ vertices, that
is, all graphs on a fixed number of vertices will have the same vertex set.
Assume $K\in\mathcal{K}(p,q)$ and let $G_{1},...,G_{l}$ and $H_{1},...,H_{l}$
be graphs such that $K=\bigoplus_{k=1}^{l}(G_{k}\otimes H_{k})$. Thus, by the
above assumption, $V(K)=\{(g_{i},h_{j})\ |\ 1\leq i\leq p,1\leq j\leq q\}$.
The next notions will be useful for Theorem \ref{thm:2char}:

\begin{definition}
Let $K$ be a (spanning) subgraph of the tensor product $G\otimes H$. Then $K$
is a \emph{cross-like subgraph} if $\{(g_{i},h_{j}),(g_{i^{\prime}%
},h_{j^{\prime}})\}\in E(K)$ implies that $\{(g_{i},h_{j^{\prime}%
}),(g_{i^{\prime}},h_{j})\}\in E(K)$ as well.
\end{definition}

\begin{definition}
Let $G$ and $H$ be graphs on vertex sets $\{g_{1},\ldots,g_{p}\}$ and
$\{h_{1},\ldots,h_{q}\}$, respectively, with $E(G)=\{\{g_{i},g_{i^{\prime}%
}\}\}$ and $E(H)=\{\{h_{j},h_{j^{\prime}}\}\}$. Let us denote the tensor
product $G\otimes H$ with $E(i,i^{\prime};j,j^{\prime})$ and call it a
\emph{tensor-elementary graph}.
\end{definition}

Using these concepts, the graphs in the set $\mathcal{K}(p,q)$ can be
characterized follows:

\begin{theorem}
\label{thm:2char} For a graph $K$, the following statements are equivalent.

(i) $K\in\mathcal{K}(p,q)$;

(ii) $K$ is a spanning, cross-like subgraph of $K_{p}\otimes K_{q}$;

(iii) $K$ is a 2-sum of tensor-elementary graphs.
\end{theorem}

\begin{proof}
\emph{(i) }$\Rightarrow(ii)$. Let $K=\bigoplus_{k=1}^{l}(G_{k}\otimes H_{k})
$. Then $V(K)=\{(g_{i},h_{j})\ |\ 1\leq i\leq p,1\leq j\leq q\}$. Consider
vertices $(g_{i},h_{j})$ and $(g_{i},h_{j^{\prime}})$ of $K$, where
$j\not =j^{\prime}$. Since $\{(g_{i},h_{j}),(g_{i},h_{j^{\prime}})\}\notin
E(G_{k}\otimes H_{k})$, $1\leq k\leq l$, we infer that $\{(g_{i},h_{j}%
),(g_{i},h_{j^{\prime}})\}\notin E(K)$. Analogously, $\{(g_{i},h_{j}%
),(g_{i^{\prime}},h_{j})\}\notin E(K)$ for any $i$ and any $j\not =j^{\prime}%
$. It follows that $K$ is a spanning subgraph of $K_{p}\otimes K_{q}$.

Suppose next that $\{(g_{i},h_{j}),(g_{i^{\prime}},h_{j^{\prime}})\}\in E(K)$.
Then $\{(g_{i},h_{j}),(g_{i^{\prime}},h_{j^{\prime}})\}\in E(G_{i}\otimes
H_{i})$ for an odd number of indices $k$, say $k=k_{1},\ldots,k_{2r+1}$,
$r\geq0$. Consequently, the edges $\{g_{i},g_{i^{\prime}}\}$ and
$\{h_{j},h_{j^{\prime}}\}$ are simultaneously present in precisely the
products $E(G_{k}\otimes H_{k})$, $k=k_{1},\ldots,k_{2r+1}$. Therefore
$\{(g_{i},h_{j^{\prime}}),(g_{i^{\prime}},h_{j})\}\in E(K)$ as well. We
conclude that $K$ is also cross-like.

\emph{(ii) }$\Rightarrow(iii)$. Let $K$ be a spanning, cross-like subgraph of
$K_{p}\otimes K_{q}$. To each pair $\{(g_{i},h_{j}),(g_{i^{\prime}%
},h_{j^{\prime}})\}$, $\{(g_{i},h_{j^{\prime}}),(g_{i^{\prime}},h_{j})\}$ of
$K$ assign the tensor-elementary graph $E(i,i^{\prime};j,j^{\prime})$. Then it
is straightforward to see that
\[
K=\bigoplus_{%
%TCIMACRO{\QATOP{\{(g_{i},h_{j}),(g_{i^{\prime}},h_{j^{\prime}})\}\in
%E(K)}{\{(g_{i},h_{j^{\prime}}),(g_{i^{\prime}},h_{j})\}\in E(K)}}%
%BeginExpansion
\genfrac{}{}{0pt}{}{\{(g_{i},h_{j}),(g_{i^{\prime}},h_{j^{\prime}})\}\in
E(K)}{\{(g_{i},h_{j^{\prime}}),(g_{i^{\prime}},h_{j})\}\in E(K)}%
%EndExpansion
}E(i,i^{\prime};j,j^{\prime}).
\]

\emph{(iii) }$\Rightarrow(i)$. This implication is obvious.
\end{proof}

\bigskip A sum modulo $2$ of tensor products is not unique. More formally,
given a tensor 2-sum graph $K$, there may be different representations of the
form $K=\bigoplus_{k=1}^{l}(G_{k}\otimes H_{k})$. This is trivially analogue
to the situation holding for density matrices, where a mixed state does not
capture all the information about the kets.

To see that a representation need not be unique it is enough to recall that
the prime factor decomposition of graph with respect to the tensor product is
not unique in the class of bipartite graphs, see~\cite{imrich00} for the
general case and~\cite{brim-05} for factorization of hypercubes. On the other
hand, the prime factor decomposition is unique for connected nonbipartite
graphs~\cite{mcke-71}, To see that this does not hold for tensor 2-sum
representations, observe first that the 2-sum is commutative and associative.
Moreover, it is not difficult to verify the distributivity law:
\begin{equation}
\label{eq:distrib}G \otimes(H_{1} \oplus H_{2}) = (G \otimes H_{1}) \oplus(G
\otimes H_{2})\,.
\end{equation}
Consider now a tensor 2-sum graph $K$ in which the first factor is fixed, that
is,
\[
K=\bigoplus_{k=1}^{l}(G\otimes H_{k})\,.
\]
Then by \eqref{eq:distrib} we can also write $K$ as
\[
K= G \otimes\left(  \oplus_{k=1}^{l} H_{k} \right)  \,.
\]
Moreover, by the commutativity and associativity of the 2-sum, the graphs
$H_{i}$ can be arbitrarily combined to get numerous different representations
of $K$.

\section{Partial transpose\label{tra}}

The Peres-Horodecki criterion for testing separability of quantum states is
based on the partial transpose of a density matrix (see, \emph{e.g.},
\cite{nc}). The criterion states that if the density matrix (or, equivalently,
the state) of a quantum mechanical system with composite dimension $pq$ is
entangled, with respect to the subsystems of dimension $p$ and $q$, then its
partial transpose is positive. For generic matrices, this operation is defined
as follows:

\begin{definition}
Let $M$ be an $n\times n$ matrix, where $n=pq$, $p,q>1$. Consider $M$ as
partitioned into $p^{2}$ blocks each of size $q\times q$. The \emph{partial
transpose} of $M$, denoted by $M^{\Gamma_{p}}$, is the matrix obtained from
$M$, by transposing independently each of its $p^{2}$ blocks. Formally,%
\[
M=\left(
\begin{array}
[c]{ccc}%
\mathcal{B}_{1,1} & \cdots & \mathcal{B}_{1,p}\\
\vdots & \ddots & \vdots\\
\mathcal{B}_{p,1} & \cdots & \mathcal{B}_{p,p}%
\end{array}
\right)  \Longrightarrow M^{\Gamma_{p}}=\left(
\begin{array}
[c]{ccc}%
\mathcal{B}_{1,1}^{T} & \cdots & \mathcal{B}_{1,p}^{T}\\
\vdots & \ddots & \vdots\\
\mathcal{B}_{p,1}^{T} & \cdots & \mathcal{B}_{p,p}^{T}%
\end{array}
\right)  ,
\]
where $\mathcal{B}_{i,j}^{T}$ denotes the transpose of the block
$\mathcal{B}_{i,j}$, for $1\leq i,j\leq p$.
\end{definition}

It is clear that we can have a partial transpose of a graph via its adjacency
matrix. The next result translates the Peres-Horodecki criterion in our
restricted setting. In a stronger way, the positivity is substituted by the
equality. This observation closely resembles the result obtained in \cite{b},
when considering normalized Laplacians. However, here we drop the constraints
of positivity and unit trace. The only property of relevance for this
criterion to hold is then symmetricity, apart from the fact that here we have
only matrices of zeros and ones.

\begin{theorem}
\label{thm:transpose} Let $K\in\mathcal{K}(p,q)$. Then $A(K)=A(K)^{\Gamma_{p}%
}$.
\end{theorem}

\begin{proof}
Let $K=\bigoplus_{k=1}^{l}(G_{k}\otimes H_{k})$. As earlier we can assume that
$V(K)=\{(g_{i},h_{j})\ |\ 1\leq i\leq p,1\leq j\leq q\}$. Also, $A(K)$ is
assumed to be constructed with respect to the lexicographic order of the
vertices of $K$: $(g_{1},h_{1}),\ldots,(g_{1},h_{q}),(g_{2},h_{1}%
),\ldots,(g_{2},h_{q}),\ldots,(g_{p},h_{q})$. To simplify the notation,
identify the vertices of $K$ in this order with the sequence $1,\ldots
,q,q+1,\ldots,2q,\ldots,pq$. Then any $i$, $1\leq i\leq pq$, can be (uniquely)
written as $i=sq+r$ for some $0\leq s\leq p-1$ and $1\leq r\leq q$. Consider
an arbitrary block $\mathcal{B}_{s_{1},s_{2}}$, $0\leq s_{1},s_{2}\leq p-1$,
of $A(K)$. Note first that by the lexicographic order, $\mathcal{B}%
_{s_{1},s_{2}}=0$ if $s_{1}=s_{2}$. Hence assume without loss of generality
$s_{1}<s_{2}$. Let the $(r_{1},r_{2})$-th entry of $\mathcal{B}_{s_{1},s_{2}}$
be equal to 1: $(\mathcal{B}_{s_{1},s_{2}})_{r_{1},r_{2}}=1$. Then
$r_{1}\not =r_{2}$. So $s_{1}q+r_{1}$ is adjacent to $s_{2}q+r_{2}$. Hence by
Theorem~\ref{thm:2char} \emph{(ii)}, $s_{2}q+r_{1}$ is adjacent to
$s_{1}q+r_{2}$. But then $(\mathcal{B}_{s_{1},s_{2}})_{r_{2},r_{1}}=1$ which
implies that $\mathcal{B}_{s_{1},s_{2}}=(\mathcal{B}_{s_{1},s_{2}})^{T}$ as claimed.
\end{proof}

\bigskip

The converse of Theorem~\ref{thm:transpose} does not hold. Consider, for
instance, the path on $4$ vertices $P_{4}$ and label its consecutive vertices
with $4,1,2,3$. Then the corresponding adjacency matrix is
\[
\left(
\begin{array}
[c]{cccc}%
0 & 0 & 1 & 1\\
0 & 0 & 1 & 0\\
1 & 1 & 0 & 0\\
1 & 0 & 0 & 0
\end{array}
\right)  ,
\]
which can be partitioned into $2\times2$ symmetric blocks. However,
$P_{4}\notin\mathcal{K}(p,q)$ since it has an odd number of edges.

\bigskip

While all separable quantum states belong to a set of PPT states (or,
\emph{Positive Partial Transpose }states), it is not immediate to construct a
general PPT state (see \cite{nc}). For graphs we have a simple method
described in the next result, where $\cup$ denotes the disjoint union of graphs.

\begin{theorem}
Let $G$ be a graph on $n$ vertices and with $m$ edges. Then the graph
\[
G\cup mK_{2} \cup(n^{2} - n - 2m)K_{1}%
\]
belongs to $\mathcal{K}(n,n)$.
\end{theorem}

\begin{proof}
Let $V(G)=\{g_{1},g_{2},\ldots,g_{n}\}$ and let $G^{\prime}$ be an isomorphic
copy of $G$ with $V(G^{\prime})=\{g_{1}^{\prime},g_{2}^{\prime},\ldots
,g_{n}^{\prime}\}$. Let $H$ be the graph with the vertex set $V(H)=\{(g_{i}%
,g_{j}^{\prime})\ |\ 1\leq i\leq p,1\leq j\leq q\}$ and the edge set
$E(H)=\{\{(g_{i},g_{i}^{\prime}),(g_{j},g_{j}^{\prime})\},\{(g_{i}%
,g_{j}^{\prime}),(g_{j},g_{i}^{\prime})\}\ |\ \{g_{i},g_{j}\}\in E(G)\}\,$.
Then it is straightforward to see that the connected components of $H$ are
$G$, $n$ copies of $K_{2}$, and the remaining $n^{2}-n-2m$ components are
$K_{1}$. In other words, $H=G\cup mK_{2}\cup(n^{2}-n-2m)K_{1}$. Moreover, $H$
is a spanning, cross-like subgraph of $K_{n}\otimes K_{n}$ so we conclude that
$H\in\mathcal{K}(n,n)$.
\end{proof}

\section{Conclusions and open problems}

In the attempt to define a minimal mathematical framework for isolating some
of the characteristic properties of quantum entanglement, we have introduced a
generalization of the tensor product of graphs. The generalization consists on
obtaining every graph by addiction modulo two, possibly with many addenda, of
tensor products of adjacency matrices. Then, we have proved a combinatorial
analogue of the Peres-Horodecki criterion, by substituting positivity with equality.

The tensor 2-sum operation gives numerous interesting issues worth of
investigation. Here is a selection of such open topics and problems.

\begin{itemize}
\item We have seen that a given graph $K$ can have (and in the most cases in
does) have different representation as a tensor 2-sum graph. Hence it is
natural to define $T_{2}(K)$ as the smallest integer $l$ (if it exists) such
that $K$ has a representation of the form $K=\bigoplus_{k=1}^{l}(G_{k}\otimes
H_{k})$. Clearly, $T_{2}(K)<\infty$ if and only if $K\in\mathcal{K}(p,q)$ for
some $p$ and $q$. The representation of $K\in\mathcal{K}(p,q)$ from
Theorem~\ref{thm:2char}~\emph{(iii)} can have arbitrarily larger number of
modulo 2 summands than $\mathrm{Kron}(K)$. Consider, for instance,
$K=K_{p}\otimes K_{q}$. Clearly, $\mathrm{Kron}(K)=1$, on the other hand the
representation of Theorem~\ref{thm:2char}~\emph{(iii)} requires $pq$ summands.
However, let $K=\bigoplus_{i=1}^{p-1}E(i,i+1;i,i+1)$. Then $T_{2}(K)=p$. Note
also that $T_{2}(K) = 1$ if and only if $K$ is not prime with respect to the
tensor product. Is there a nice characterization of graphs $K$ with $T_{2}%
(K)$? More generally, it would be nice to have a classification of graphs in
terms of the minimum number of summands required for their constructions as a
sum of tensor products (that is, in terms of $T_{2}$).

\item Theorem~\ref{thm:2char} gives two necessary and sufficient conditions
for a graph to belong to $\mathcal{K}(p,q)$. However, these conditions are not
efficient, so it remains to determine the computational complexity of the
following decision problem:

\begin{itemize}
\item \textbf{Given:} A graph $G$ on $n=pq$ vertices.

\item \textbf{Task:} Is $G\in\mathcal{K}(p,q)$?
\end{itemize}

We feel that recent investigations of the so-called approximate graph
products~\cite{heim-09} might be useful in solving this problem. In this
respect we add that the unique prime factorization of nonbipartite connected
graphs can be found in polynomial time~\cite{im-98}.

\item Suppose $K$ is a tensor 2-sum graph with a representation $K=\bigoplus
_{k=1}^{l}(G_{k}\otimes H_{k})$. Then the only condition we posed on the
graphs $G_{i}$ and $H_{i}$ is that each has at least one edge. One might want
to be also more restrictive by imposing that all $G_{i}$'s and $H_{i}$'s must
be connected. What can be said of such restricted representations?
\end{itemize}

\bigskip

\emph{Acknowledgements. }Sandi Klav\v{z}ar is supported in part by the
Ministry of Science of Slovenia under the grant P1-0297. The author is also
with the Institute of Mathematics, Physics and Mechanics, Jadranska 19, 1000
Ljubljana. Simone Severini is a Newton International Fellow.

\end{document}